\documentclass[12pt]{article}
\usepackage[colorlinks,
            linkcolor=blue,
            anchorcolor=green,
            citecolor=blue
            ]{hyperref}
\usepackage{mathrsfs,amssymb,amsfonts}
\usepackage{amsmath,amssymb,latexsym,color}
\usepackage[mathscr]{eucal}
\usepackage[numbers,sort&compress]{natbib}
\usepackage{graphicx}
\usepackage{lipsum}
\usepackage{enumerate}
\usepackage[utf8]{inputenc}
\def\thefootnote{\fnsymbol{footnote}}
\newtheorem{thm}{Theorem}[section]

\newtheorem{prop}[thm]{Proposition}
\newtheorem{lemma}[thm]{Lemma}
\newtheorem{cor}[thm]{Corollary}
\newtheorem{example}[thm]{Example}

\newtheorem{problem}[thm]{Problem}
\newtheorem{ques}[thm]{Question}
\newtheorem{remark}[thm]{Remark}
\newtheorem{Notation}[thm]{Notation}
\newtheorem{obs}[thm]{Observation}

\newcommand{\proof}{{\it Proof.\quad}}
\newcommand{\qed}{\hfill\Box\medskip}
\renewcommand{\thefootnote}{\arabic{footnote}}

\begin{document}
\renewcommand{\abovewithdelims}[2]{
\genfrac{[}{]}{0pt}{}{#1}{#2}}
%%%%%%%%%%%%%%%%%%%%%%%%%%%%%%%%%%%%%%%%%%%%%%%%%%%%%%%%%%%%%%%%%%%%%%%%%%%%%%%%%%%%%%%%
%%%%%%%%%%%%%%%%%%%%%%%%%%%%%%%%%%%%%%%%%%%%%%%%%%%%%%%%%%%%%%%%%%%%%%%%%%%%%%%%%%%%%%%%

\title{\bf On finite groups whose coprime graph is a divisor graph
}

\author{{\bf Xuanlong Ma$^{{\rm 1}}$, Liangliang Zhai$^{{\rm 1}}$, Nan Gao$^{{\rm 1}}$, Junyao Pan$^{{\rm 2,}}$\footnote{Corresponding author}}\\
{\footnotesize 1. School of Science, Xi'an Petroleum University, Xi'an 710065, China}\\
{\footnotesize 2. School of Science, Guangxi University of Science and Technology, Liuzhou 545005, China}}

 \date{}
 \maketitle
\newcommand\blfootnote[1]{%
\begingroup
\renewcommand\thefootnote{}\footnote{#1}%
\addtocounter{footnote}{-1}%
\endgroup
}
\begin{abstract}
In this paper, we first characterize which generalized lexicographic products are divisor graphs.
As applications, we show that power graphs, reduced power graphs and order graphs are all divisor graphs, which also implies the main result in
[Power graph of a finite group is always divisor graph, Asian-European Journal of Mathematics 16 (2023)]. We then show that, the coprime graph of a group is a generalized
lexicographic product, and characterize which coprime graphs are divisor graphs. Finally, we classify the finite groups $G$ having at most four prime divisors, whose coprime graphs are divisor graphs, and we also classify the finite groups $G$ whose coprime graphs are divisor graphs, if $G$ is a nilpotent group, a dihedral group, a generalized quaternion group, a symmetric group, an alternating group, a direct product of two non-trivial groups, and a sporadic simple group.


\end{abstract}

\blfootnote{{\em E-mail addresses:}
xuanlma@xsyu.edu.cn (Xuanlong Ma), zhailiang111@126.com (Liangliang Zhai), gaonankk@163.com (Nan Gao), Junyao$_{-}$Pan@126.com (Junyao Pan)
}

\section{Introduction}

The study of graph representations of an algebraic structure is
a popular and interesting research topic in the field of algebraic graph theory.
For example, a well-known graph representation for group is
the Cayley graph, which has a very long history.
Moreover, graphs from algebraic structures have been actively investigated in the literature, since they have valuable applications (cf. \cite{K}).

Given a group $G$, one can define various graphs on $G$.
Among a number of graphs defined on the
group $G$ reflecting some algebraic properties of $G$,
apart from Cayley graph,
two which have been widely studied are the following:
\begin{itemize}
  \item The {\em commuting graph} of $G$ is the undirected graph with vertex set $G$, and two distinct vertices are adjacent if they commute,  
  which was first introduced by Brauer and Fowler \cite{B55};
  \item The {\em power graph} of $G$, denoted by $\mathcal{P}(G)$, is the undirected graph with vertex set $G$, and two distinct vertices are adjacent if one is a power of the other.
\end{itemize}
In 2000, Kelarev and Quinn \cite{KQ00} first introduced the concept of a power graph defined on a group, which is a directed graph.
In 2009, Chakrabarty {\em et al.} \cite{CGS} introduced the concept of an undirected power graph. Afterwards, the term ``power graph" always means  an undirected power graph (cf. \cite{AKC}).

Because the order of an element is one of the most basic and important concepts in group theory, one can define a graph over a group by element order. For a group $G$,
the {\em order graph} $\mathcal{S}(G)$ of $G$ is an undirected graph with vertex set $G$, and two distinct vertices $x,y$ are adjacent if $o(x)|o(y)$ or $o(y)|o(x)$, where $o(x)$ and $o(y)$ are the orders of $x$ and $y$, respectively.
In 2017, Hamzeh and Ashrafi \cite{HA17} first introduced the order graph and called this graph as the main supergraph of $\mathcal{P}(G)$. Clearly, if two distinct vertices $a$ and $b$ are adjacent in $\mathcal{P}(G)$, then $a$ and $b$ are also adjacent in $\mathcal{S}(G)$. It follows that $\mathcal{P}(G)$ is a spanning subgraph of $\mathcal{S}(G)$.

In \cite{Ma13}, the authors introduced the concept of a coprime graph of a group. For a group $G$, the {\em coprime graph} of $G$, denoted by $\Gamma(G)$, is the undirected graph with vertex set $G$ and two distinct vertices $x,y$ are adjacent if $o(x)$ and $o(y)$ are relatively prime, that is, $(o(x),o(y))=1$.
The authors in \cite{Ma13} explored how the
properties of graphs can effect on the properties of groups.
Dorbidi \cite{Do} showed that, for every finite group $G$, the clique number of $\Gamma(G)$ is equal to the chromatic number of $\Gamma(G)$, and classified the groups whose coprime graph is complete $r$-partite or planar. Selvakumar and Subajini \cite{Se} classified the finite groups whose coprime graph has genus one.
Hamm and Way \cite{Ha} obtained the exact value of the independence number of the coprime graph of a dihedral group and studied the finite groups whose  coprime graph is perfect.
Alraqad {\em et al.} \cite{Al} obtained all finite groups whose coprime graph has precisely $3$ end-vertices.

Every graph considered in this paper is simple, without loops and multiple edges. Let $\Gamma$ be a graph. Denote by $V(\Gamma)$ and $E(\Gamma)$ the vertex set and edge set of $\Gamma$, respectively.
An {\em orientation} $D$ of an undirected graph $\Gamma$ is an assignment of exactly one direction to each edge of $\Gamma$.
Let $D$ be an orientation of $\Gamma$.
A {\em transmitter} of $D$ is a vertex having indegree $0$ and a {\em receiver} of $D$ is a vertex having outdegree $0$. A vertex $b$ of $D$ is {\em transitive} provided that its outdegree and indegree are both greater than zero, and $(a,b),(b,c)\in E(D)$ must imply $(a,c)\in E(D)$, where $a,c\in V(D)$ and $(a,b)\in E(D)$ means that there is a  directed edge from $a$ to $b$ in $D$.
A {\em graph labeling} of a graph is an assignment of integers to the vertices or edges or both, subject to certain conditions have been motivated by practical problems. Labeled graphs serve useful mathematical models for a broad range of applications, such as, satellite communication through \cite{P} and life sciences \cite{Fo}.

Singh and Santhosh \cite{sin} first introduced the concept of a
divisor graph.
A {\em divisor graph} $\Gamma$ is a simple graph with $V(\Gamma)$ as a subset of the set of positive integers, and any two distinct $x,y\in V(\Gamma)$, $\{x,y\}\in E(\Gamma)$ if and only if either $x\mid y$ or $y\mid x$.
A graph isomorphic to a divisor graph
is also called a {\em divisor graph}.
That is to say, if $\Gamma$ is a divisor graph, then there exists an injection $\alpha$ from $V(\Gamma)$ to the set of positive integers such that $\Gamma$ is isomorphic to the divisor graph with vertex set $\alpha(V(\Gamma))$, where $\alpha$ is called a {\em divisor labeling} of $\Gamma$. In 2001, Chartrand {\em et al.} \cite{Char} studied the following question:
\begin{ques}\label{ques1}
Which graphs are divisor graphs?
\end{ques}
Specifically, they showed that, a graph $\Gamma$ is a divisor graph if and only if there exists an orientation $D$ of $\Gamma$ such that every vertex of $D$ is a transmitter, a
receiver, or a transitive vertex. Afterward, Frayer \cite{Fr} also studied Question~\ref{ques1} and showed that,
for a divisor graph $\Gamma$ with
a transitive vertex, $\Gamma \times \Delta$ is a divisor graph if and only if $E(\Delta)=\emptyset$.
In 2006, Vinh proved that a graph $\Gamma$ is a divisor graph if and only if $\Gamma$ has a transitive orientation (see Theorem~\ref{dg-o}).
In 2012, Al-Addasi {\em et al.} \cite{Ala} characterized the graphs $\Gamma$ and $\Delta$ for which the Cartesian product of $\Gamma$ and $\Delta$ is a divisor graph. They also show that divisor graphs form a proper subclass of perfect graphs.
Number theoretic graphs are one of the emerging fields in graph theory.
Recently, Ravi and Desikan \cite{RD23} published a brief survey on divisor graphs and divisor function graphs.

The zero-divisor graph is very famous among all graphs associated with a ring, which was first introduced by Beck \cite{Be}. Osba and Alkam \cite{Os} studied this question: Which zero-divisor graphs are divisor graphs?
Motivated by Question~\ref{ques1}, in this paper, we will study the following question:

\begin{ques}\label{ques2}
Which coprime graphs are divisor graphs?
\end{ques}

In this paper, we first characterize which generalized lexicographic products are divisor graphs.
As applications, we show that every of power graph, reduced power graph and order graph is a divisor graph. We then show that, the coprime graph of a group is a generalized
lexicographic product, and characterize which
coprime graphs are divisor graphs.
Finally, we classify the finite groups $G$ having at most four prime divisors, whose coprime graphs are divisor graphs.
We also classify the finite groups $G$ whose coprime graphs are divisor graphs, if $G$ is a nilpotent group, a dihedral group, a generalized quaternion group, a symmetric group, an alternating group, a direct product of two non-trivial groups, and a sporadic simple group.
These results partly answer Question~\ref{ques2}.

\section{Preliminaries}\label{}
In this section, we briefly introduce some notation, terminology, and basic results that will be used in the sequel.

Denote by $K_n$ the complete graph of order $n$.
Let $\Gamma$ and $\Delta$ be two graphs. If $V(\Delta)\subseteq V(\Gamma)$ and for distinct $a,b\in V(\Delta)$, it must be that $\{a,b\}\in E(\Gamma)$ implies $\{a,b\}\in E(\Delta)$, then $\Delta$ is called an {\em induced subgraph} of $\Gamma$.
In particular, if $K$ is a subset of $V(\Gamma)$, then the induced subgraph by $K$ is denoted by $\Gamma[K]$.
If $\Gamma_1$ and $\Gamma_2$ are two graphs with $V(\Gamma_1)\cap V(\Gamma_2)=\emptyset$, then the {\em union} of $\Gamma_1$ and $\Gamma_2$, denoted by $\Gamma_1 \cup \Gamma_2$, is the graph with vertex set $V(\Gamma_1)\cup V(\Gamma_2)$ and edge set $E(\Gamma_1)\cup E(\Gamma_2)$.
For two graphs $\Gamma_1$ and $\Gamma_2$ with $V(\Gamma_1)\cap V(\Gamma_2)=\emptyset$, the {\em sum} of $\Gamma_1$ and $\Gamma_2$, denoted by $\Gamma_1 \vee  \Gamma_2$, is the graph with vertex set $V(\Gamma_1)\cup V(\Gamma_2)$, and its edge set is the union of $E(\Gamma_1)$, $ E(\Gamma_2)$ and the set of
all edges between vertices from the two different graphs $\Gamma_1$ and $\Gamma_2$.

The next three results are derived from \cite{Char}.

\begin{lemma}\label{indsub}
Every induced subgraph of a divisor graph is also a divisor graph.
\end{lemma}

\begin{lemma}\label{complete}
Every complete graph is a divisor graph.
\end{lemma}

\begin{lemma}\label{liantu}
If both $\Gamma_1$ and $\Gamma_2$ are divisor graphs, 
then both 
$\Gamma_1 \vee \Gamma_2$ and $\Gamma_1 \cup \Gamma_2$ are also divisor graphs. 
\end{lemma}

The following result is derived from \cite{Char} and also
appears in \cite{Fr}.

\begin{lemma}{\rm (\cite[Theorem 3.1 ]{Fr})}\label{bip}
Every bipartite graph is a divisor graph.
\end{lemma}

Vinh \cite{Vi} showed that divisor graphs have arbitrary order and size. In particular, he gave a simple proof for the following theorem. Indeed, the following depiction of divisor graphs is due
to Chartran, Muntean, Saenpholpant and Zhang \cite{Char}.

\begin{thm}{\rm (\cite[Theorem 2]{Vi})}\label{dg-o}
A graph $\Gamma$ is a divisor graph if and only if there exists an orientation $\mathcal{D}$ on $\Gamma$ such that if $(x, y), (y, z)$ are directed edges of $\mathcal{D}$ then so is $(x, z)$, that is, $\Gamma$ is a transitive orientation of $\Gamma$.
\end{thm}

Every group considered in our paper is finite.
We always use $e$ to denote the identity
of the group under consideration.
A finite group is called a {\em CP-group} provided that its
every element is of prime power order (cf. \cite{DW}).
For example, the alternating group on $5$ letters is a CP-group.
Certainly, every $p$-group is also a CP-group for some prime $p$.
Recall that Delgado and Wu characterized
the finite CP-groups (see \cite[Theorem 4]{DW}).
Given a finite group $G$, denote by $\pi_e(G)$ and $\pi(G)$
the set of orders of all non-identity elements of $G$ and the set of all prime divisors of $|G|$, respectively.
For $a\in G$, if the situation is unambiguous, then we denote $\pi(\langle a\rangle)$ simply by $\pi(a)$.
The {\em symmetric group} of degree $n$ and the {\em alternating group} of degree $n$ are denoted by 
$\mathbf{S}_n$ and $\mathbf{A}_n$, respectively. 

By the definition of a coprime graph, if $H$ is a subgroup of a group $G$, then $\Gamma(H)$ is an induced subgraph of $\Gamma(G)$. So,
the following result holds.

\begin{obs}\label{obs-1}
Give a group $G$, $\Gamma(G)$ is a divisor graph if and only if for every subgroup $H$ of $G$, $\Gamma(H)$ is a divisor graph.
\end{obs}

\section{Generalized Lexicographic Products}\label{}

In this section, we give a characterization for which generalized lexicographic products are divisor graphs. As applications, we show that every of power graph, reduced power graph and order graph is a divisor graph.

We first give the definition of a generalized lexicographic
product, which was first introduced by Sabidussi \cite{Sa61}.
Let $\mathcal{H}$ be a graph. For every $v\in V(\mathcal{H})$, we use $\mathcal{F}_v$ to denote a graph indexed by the vertex $v$. Let $\mathbb{F}$ be a family of graphs as follows:
$$\mathbb{F}=\{\mathcal{F}_v: v\in V(\mathcal{H})\}.$$
The {\em generalized lexicographic product} of $\mathcal{H}$ and $\mathbb{F}$, denoted by $\mathcal{H}[\mathbb{F}]$, is an undirected graph with vertex set $\{(v,w): v\in V(\mathcal{H}), w\in V(\mathcal{F}_v)\}$ and edge set
$$\{\{(v_1,w_1),(v_2,w_2)\}:\{v_1,v_2\}\in E(\mathcal{H}),
\text{ or $v_1=v_2$ and } \{w_1,w_2\}\in E(\mathcal{F}_{v_1})\}.$$

\begin{thm}\label{glp}
$\mathcal{H}[\mathbb{F}]$ is a divisor graph if and only if $\mathcal{H}$ is a divisor graph and for every $v\in V(\mathcal{H})$, $\mathcal{F}_v$ is also a divisor graph.
\end{thm}
\proof
Suppose first that $\mathcal{H}$ is a divisor graph and $\mathcal{F}_v$ is also a divisor graph for every $v\in V(\mathcal{H})$.
By Theorem~\ref{dg-o}, let $\mathcal{D}'$ be a transitive orientation of $\mathcal{H}$, and let $\mathcal{D}_v$ be a transitive orientation of $\mathcal{F}_v$ for every $v\in V(\mathcal{H})$.
For $\{v_1,v_2\}\in E(\mathcal{H})$, $w_1\in V(\mathcal{F}_{v_1})$ and $w_2\in V(\mathcal{F}_{v_2})$, let the direction of
$\{(v_1,w_1),(v_2,w_2)\}$ be the direction of $\{v_1,v_2\}$ in $\mathcal{D}'$. Also, for $\{w_1,w_2\}\in E(\mathcal{F}_{v_1})$,
let the direction of
$\{(v_1,w_1),(v_1,w_2)\}$ be the direction of $\{w_1,w_2\}$ in $\mathcal{D}_{v_1}$. Then we obtain a orientation of $\mathcal{H}[\mathbb{F}]$, say $\mathcal{D}$. It is clear that $\mathcal{D}$ is transitive, and so $\mathcal{H}[\mathbb{F}]$ is a divisor graph by Theorem~\ref{dg-o}.

For the converse, let $\mathcal{H}[\mathbb{F}]$ be a divisor graph and let $\mathcal{D}$ denote its a transitive orientation. Then it is easy to see that the subgraph of $\mathcal{H}[\mathbb{F}]$ by
$\{(v,w): w\in V(\mathcal{F}_{v})\}$ is isomorphic to $\mathcal{F}_{v}$, it follows from Lemma~\ref{indsub} that $\mathcal{F}_{v}$ is a divisor graph for every $v\in V(\mathcal{H})$. Also, for every
$v\in V(\mathcal{H})$,
fix a vertex, say $w_v$, in $V(\mathcal{F}_{v})$. Then
the subgraph of $\mathcal{H}[\mathbb{F}]$ induced by
$\{(v,w_v): v\in V(\mathcal{H})\}$ is isomorphic to $\mathcal{H}$, which implies that $\mathcal{H}$ is a divisor graph.
$\qed$

Given a finite group $G$, define a relation $\equiv$ on $G$ as follows:
$$
x\equiv y \Leftrightarrow \langle x\rangle=\langle y\rangle.
$$
Clearly, $\equiv$ is an equivalence relation over $G$. Denote by $\overline{x}$ the $\equiv$-class containing $x\in G$.
Let $\overline{G}=\{\overline{x}:x\in G\}$.
Now we define the undirected graph $\mathcal{L}_G$ with vertex set $\overline{G}$, and for distinct $\overline{x},\overline{y}\in \overline{G}$, $\{\overline{x},\overline{y}\}\in E(\mathcal{L}_G)$ if and only if either $\langle x\rangle \subsetneq \langle y\rangle$ or $\langle y\rangle \subsetneq \langle x\rangle$.
For any $\{\overline{x},\overline{y}\}\in E(\mathcal{L}_G)$ with $\langle x\rangle \subsetneq \langle y\rangle$,
then from $\overline{y}$ to $\overline{x}$, we give the direction of $\{\overline{x},\overline{y}\}$. Thus, we can get an orientation of $\mathcal{L}_G$, which is transitive. So, $\mathcal{L}_G$ is a divisor graph by Theorem~\ref{dg-o}.

In \cite{FMW}, the authors showed that every power graph is isomorphic to the generalized lexicographic product $\mathcal{L}_G[\mathbb{K}_G]$, where every graph in $\mathbb{K}_G$ is complete (see \cite[Theorem 2.14]{FMW}).
Note that by Lemma~\ref{complete}, every complete graph is a divisor graph. Thus, as a corollary of Theorem~\ref{glp}, we have the following result, which is also the main result of \cite{Ta} (see \cite[Theorem 4]{Ta}). Actually, in \cite{FMW}, the authors have given a transitive orientation of $\mathcal{P}(G)$ (see \cite[Theorem 2.3]{FMW}).

\begin{cor}
The power graph of a finite group is a divisor graph.
\end{cor}

In order to avoid the complexity of edges in the power graphs, Rajkumar and Anitha \cite{RP1} introduced the {\em reduced power graph} $\mathcal{P}_{R}(G)$ of a finite group $G$, which is an undirected graph with vertex
set $G$, and two distinct vertices $x,y$  are adjacent if
$\langle x\rangle \subset \langle y\rangle$ or $\langle y\rangle \subset \langle x\rangle$. In fact, $\mathcal{P}_{R}(G)$ is the subgraph of $\mathcal{P}(G)$ obtained by deleting the edges $\{x,y\}$ satisfying  $\langle x\rangle=\langle y\rangle$. It is similar to power graphs,
one can easily obtain that $\mathcal{P}_{R}(G)$ is isomorphic to the generalized lexicographic product $\mathcal{L}_G[\mathbb{I}_G]$, where every graph in $\mathbb{I}_G$ is empty. Clearly, any empty is also a divisor graph. Thus, by Theorem~\ref{glp}, we have the following result, which is also a result in \cite{TS}. 

\begin{cor}
The reduced power graph of a finite group is a divisor graph.
\end{cor}

For $x,y \in G$, define
$x\thickapprox y$ if $o(x)=o(y)$. It is readily seen that $\thickapprox$ is an equivalence relation on $G$. We denote by $\widehat{x}$
the $\thickapprox$-class containing $x$.
Clearly, $\thickapprox$ is a clique in $\mathcal{S}_G$.
Write $\widehat{G}=\{\widehat{x}: x\in G\}$.
Now we define an undirected graph $\mathcal{J}_G$ with vertex set $\widehat{G}$, and for distinct $\widehat{x},\widehat{y}\in \widehat{G}$, $\{\widehat{x},\widehat{y}\}\in E(\mathcal{J}_G)$ if and only if either $o(x)\mid o(y)$ or $o(y)\mid o(x)$.
For any $\{\widehat{x},\widehat{y}\}\in E(\mathcal{J}_G)$ with $o(x)\mid o(y)$,
then from $\widehat{x}$ to $\widehat{y}$, we give the direction of $\{\widehat{x},\widehat{y}\}$. Clearly, we obtain a transitive orientation of $\mathcal{J}_G$. So, by Theorem~\ref{dg-o}$, \mathcal{J}_G$ is a divisor graph.
Ma and Su showed that $\mathcal{S}(G)$ is isomorphic to the generalized
lexicographic product $\mathcal{J}_G[\mathbb{K}_G]$, where every graph in $\mathbb{K}_G$ is complete (see \cite[Theorem 2.1]{MS}). Thus, the following holds.

\begin{cor}
The order graph of a finite group is a divisor graph.
\end{cor}

\section{Coprime Graphs}\label{}

In this section, we first prove that the coprime graph of a group is a generalized lexicographic product (see Proposition~\ref{str}). Then, we will characterize which coprime graphs are divisor graphs (see Theorem~\ref{cp-div}).

In the following,
we always use $G$ to denote a finite group and use
$e$ to denote the identity
of $G$.

\begin{lemma}\label{co-dig}
$\Gamma(G)$ is a divisor graph if and only if
the induced subgraph $\Gamma(G)[G\setminus\{e\}]$ of $\Gamma(G)$ by $G\setminus\{e\}$ is a divisor graph.
\end{lemma}
\proof
Note that $e$ is adjacent to every vertex of  $G\setminus\{e\}$ in $\Gamma(G)$. Therefore, $\Gamma(G)\cong K_1\vee \Gamma(G)[G\setminus\{e\}]$. Now the desired result follows from
Lemmas~\ref{indsub} and \ref{liantu}.
$\qed$

Now, we define a binary relation $\circeq$ on $G\setminus\{e\}$ by the following rule:
$$
a \circeq b \Leftrightarrow \pi(a)=\pi(b),~~a,b\in G\setminus\{e\}.
$$
It is readily seen that $\circeq$ is an equivalence relation on $G\setminus\{e\}$.
For $a\in G\setminus\{e\}$, denote by $[a]$ the $\circeq$-class having $a$.
It is clear that $[a]$ is an independent set of $\Gamma(G)$.
Write $$\widetilde{G}=\{[a]: a\in G\setminus\{e\}\}.$$
Now we define a new undirected simple graph $\mathcal{O}_G$ on $\widetilde{G}$, which has vertex set $\widetilde{G}$ and edge set
$$
\{\{[a],[b]\}: [a],[b]\in \widetilde{G}\text{ and } \pi(a)\cap \pi(b)=\emptyset\}.
$$
For any $[a]\in \widetilde{G}$, let $\mathcal{F}_{[a]}$ be the induced subgraph of $\Gamma(G)$ by $[a]$. Clearly, we have $\mathcal{F}_{[a]}\cong \overline{K_{|[a]|}}$, which is the empty graph of order $|[a]|$. Let
$$\mathbb{L}_G=\{\mathcal{F}_{[a]}: [a]\in V(\mathcal{O}_G)\}.$$
Then we have the following result, which gives a characterization of the structure of a coprime graph.

\begin{prop}\label{str}
$\Gamma(G)[G\setminus\{e\}]\cong\mathcal{O}_G[\mathbb{L}_G]$.
In particular, $\Gamma(G)$ is also a generalized lexicographic product, and $\Gamma(G)\cong K_1 \vee\mathcal{O}_G[\mathbb{L}_G]$.
\end{prop}

We use the following example to illustrate Proposition~\ref{str}. 

\begin{example}
Consider the symmetric group $\mathbf{S}_3$. 
Then 
$$\mathbf{S}_3=\{(1)\}\cup \{(1,2),(1,3),(2,3)\}\cup \{(1,2,3),(1,3,2)\}$$
and $[(1,2)]=\{(1,2),(1,3),(2,3)\},[(1,2,3)]=\{(1,2,3),(1,3,2)\}$.
It follows from Proposition~\ref{str} that $\Gamma(\mathbf{S}_3)\cong K_1 \vee K_2[\{\overline{K_2},\overline{K_3}\}]$, which is shown in Figure~\ref{nt}.
\begin{figure}[hptb]
  \centering
  \includegraphics[width=5.5cm]{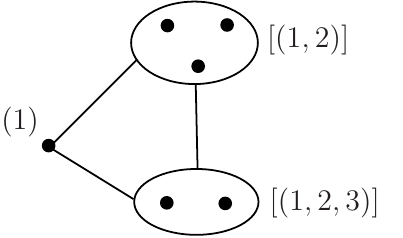}\\
  \caption{$\Gamma(\mathbf{S}_3)\cong K_1 \vee K_2[\{\overline{K_2},\overline{K_3}\}]$}\label{nt}
\end{figure}
\end{example}

We now define a new graph $\mathcal{N}_G$ with respect to $G$ as follows:
$$
V(\mathcal{N}_G)=\{\prod_{p\in \pi(a)}p: [a]\in \widetilde{G}\},
$$
and
$$
E(\mathcal{N}_G)=\{\{x,y\}: x\ne y\in V(\mathcal{N}_G)\text{ and } (x,y)=1\}.
$$
It is clear that
\begin{equation}\label{eq-01}
\mathcal{N}_G\cong \mathcal{O}_G.
\end{equation}

Note that every graph in $\mathbb{L}_G$ is a divisor graph. Therefore, by Theorem~\ref{glp}, Lemma~\ref{co-dig} and \eqref{eq-01}, we have the following result.

\begin{thm}\label{cp-div}
$\Gamma(G)$ is a divisor graph, if and only if $\mathcal{O}_G$ is a divisor graph, if and only if $\mathcal{N}_G$ is a divisor graph.
\end{thm}

Clearly, if $G$ is a CP-group, then $\mathcal{N}_G$ is a complete graph, and thus $\mathcal{N}_G$ is a divisor graph by Lemma~\ref{complete}. Now
as a direct consequence of Theorem~\ref{cp-div}, we have the following corollary.

\begin{cor}\label{cp-g}
$\Gamma(G)$ is a divisor graph for every CP-group $G$.
\end{cor}

\section{Groups $G$ with $|\pi(G)|\le 4$}\label{}
In this section, we characterize the groups $G$ with $|\pi(G)|\le 4$, whose coprime graph is a divisor graph (see Theorems~\ref{2-p-ds}, \ref{3-p-ds} and \ref{4-p-ds}). We also classify all nilpotent groups whose coprime graph is a divisor graph (see Proposition~\ref{nil-g}).

We begin with the following lemma.

\begin{figure}[hptb]
  \centering
  \includegraphics[width=7cm]{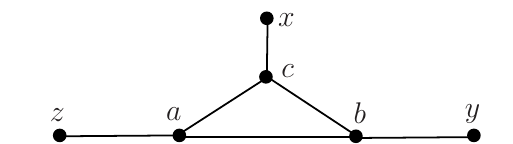}\\
  \caption{A graph which is not a divisor graph}\label{t-1}
\end{figure}

\begin{lemma}\label{teshu}
If $\{p,q,r,pq,pr,qr\}\subseteq\pi_e(G)$, where $p,q,r$ are pairwise distinct primes, then $\Gamma(G)$ is not a divisor graph.
\end{lemma}
\proof
Take $a,b,c,x,y,z\in G$ with $$o(a)=p,o(b)=q,o(c)=r,o(x)=pq,o(y)=pr,o(z)=qr.$$
Note that the subgraph of $\Gamma(G)$ induced by $\{a,b,c,x,y,z\}$
is shown in Figure~\ref{t-1}. Figure~\ref{t-1} was given as an example of a non-divisor graph in \cite{Char}.
It follows from Lemma~\ref{indsub} that $\Gamma(G)$ is not a divisor graph.
$\qed$

\begin{cor}\label{3-n-pr}
If $G$ has an element of order $pqr$ where $p,q,r$ are pairwise distinct primes, then $\Gamma(G)$ is not a divisor graph.
\end{cor}

\begin{thm}\label{2-p-ds}
If $|G|$ has at most two prime divisors, then $\Gamma(G)$ is a divisor graph.
\end{thm}
\proof
If $G$ is a $p$-group, then $\Gamma(G)\cong K_{1,|G|-1}$, and so $\Gamma(G)$ is a divisor graph by Lemma~\ref{bip}. Now suppose that $\pi(G)=\{p,q\}$ where $p,q$ are distinct primes. Then it is easy to see that
\begin{center}
$\widetilde{G}=\{[a],[b]\}$ or $\{[a],[b],[c]\}$,
\end{center}
where $\pi(a)=\{p\}$, $\pi(b)=\{q\}$ and $\pi(c)=\{p,q\}$. It follows that $\mathcal{O}_G$ is isomorphic to $K_2$ or $K_1\cup K_2$. By Lemmas~\ref{complete} and \ref{liantu}, we have that $\mathcal{O}_G$ is a divisor graph. We then derive from Theorem~\ref{cp-div} that $\Gamma(G)$ is a divisor graph, as desired.
$\qed$

Clearly, for a nilpotent group $G$, if $\pi(G)=\{p,q,r\}$,
then $G$ must have an element of order $pqr$. Now combining Corollary~\ref{3-n-pr} and Theorem~\ref{2-p-ds}, we have the following result, which classifies the nilpotent groups whose coprime graph is a divisor graph.

\begin{cor}\label{nil-g}
Given a nilpotent group $G$, $\Gamma(G)$ is a divisor graph if and only if $G$ is isomorphic to either a $p$-group or $P\times Q$, where $P$ and $Q$ are respectively a $p$-group and a $q$-group for distinct primes $p,q$.
\end{cor}

\begin{figure}[hptb]
  \centering
  \includegraphics[width=11cm]{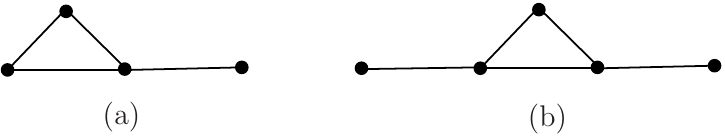}\\
  \caption{Two possibilities for $\mathcal{O}_G$}\label{t-2}
\end{figure}

\begin{figure}[hptb]
  \centering
  \includegraphics[width=11cm]{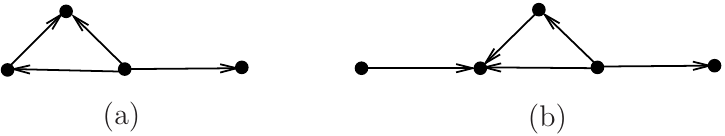}\\
  \caption{Two possible transitive orientations for $\mathcal{O}_G$}\label{t-3}
\end{figure}

\begin{thm}\label{3-p-ds}
Suppose that $\pi(G)=\{p,q,r\}$.
Then $\Gamma(G)$ is a divisor graph if and only if $\{pq,pr,qr\}\nsubseteq\pi_e(G)$.
\end{thm}
\proof
The necessity follows trivially from Lemma~\ref{teshu}. We next prove the
sufficiency. Suppose that $\{pq,pr,qr\}\nsubseteq\pi_e(G)$.
Then $G$ has no elements of order $pqr$. Furthermore, we have
$$
\{[a],[b],[c]\}\subseteq \widetilde{G}\subsetneq \{[a],[b],[c],[x],[y],[z]\},
$$
where $\pi(a)=\{p\}$, $\pi(b)=\{q\}$, $\pi(c)=\{r\}$,
$\pi(x)=\{p,q\}$, $\pi(y)=\{p,r\}$, and $\pi(z)=\{qr\}$.
It follows that $\mathcal{O}_G$ is isomorphic to one graph in Figure~\ref{t-2}. We then can obtain transitive orientations for the two graphs in Figure~\ref{t-2}, which is shown in Figure~\ref{t-3}. As a result, Theorem~\ref{dg-o} implies that $\mathcal{O}_G$ is a divisor graph, and so from Theorem~\ref{cp-div}, it follows that $\Gamma(G)$ is a divisor graph, as desired.
$\qed$

\begin{figure}[hptb]
  \centering
  \includegraphics[width=9cm]{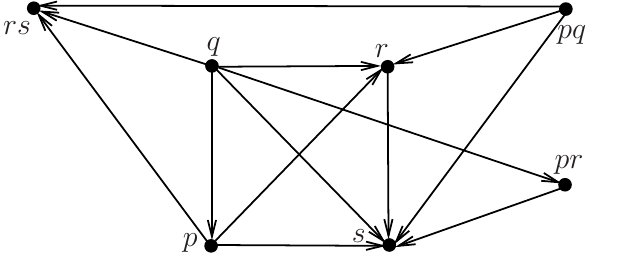}\\
  \caption{A transitive orientation for $\mathcal{O}_G$}\label{t-4}
\end{figure}

\begin{thm}\label{4-p-ds}
Suppose that $\pi(G)=\{p,q,r,s\}$.
Then $\Gamma(G)$ is a divisor graph if and only if
\begin{equation}\label{eq-new2}
\{p,q,r,s\}\subseteq V(\mathcal{N}_G)\subseteq \{p,q,r,s,pq,pr,rs\},
\end{equation}
\end{thm}
\proof
We first prove the sufficiency.
Suppose that $V(\mathcal{N}_G)\subseteq \{p,q,r,s,pq,rs,pr\}$.
We then show that if $V(\mathcal{N}_G)=\{p,q,r,s,pq,rs,pr\}$, then $\Gamma(G)$ is a divisor. In fact, in this case, $\mathcal{N}_G$ is
the underlying graph of the graph as shown in Figure~\ref{t-4}.
Clearly, Figure~\ref{t-4} is a transitive orientation of $\mathcal{N}_G$, it follows from Theorem~\ref{cp-div} that $\Gamma(G)$ is a divisor graph, as desired. Now if $V(\mathcal{N}_G)\subseteq \{p,q,r,s,pq,rs,pr\}$, then $\mathcal{N}_G$ must be an induced subgraph of the underlying graph of the graph as shown in Figure~\ref{t-4}, and so $\mathcal{N}_G$ is also a divisor graph by Lemma~\ref{indsub}. Then Theorem~\ref{cp-div} implies that $\Gamma(G)$ is a divisor graph, as desired.

Conversely, suppose that $\Gamma(G)$ is a divisor graph.
It follows from Lemma~\ref{teshu} and Corollary~\ref{3-n-pr} that
the product of pairwise distinct primes must not belong to $V(\mathcal{N}_G)$ and $\{xy,xz,yz\}\subsetneq V(\mathcal{N}_G)$ for pairwise distinct three primes $x,y,z$. As a result, without loss of generality, we may assume that either
\begin{equation}\label{eq-02}
\{p,q,r,s\}\subseteq V(\mathcal{N}_G)\subseteq \{p,q,r,s,pq,pr,rs,qs\}
\end{equation}
or
\begin{equation}\label{eq-0n2}
\{p,q,r,s\}\subseteq V(\mathcal{N}_G)\subseteq \{p,q,r,s,pq,pr,ps\}.
\end{equation}
Notice that $\mathcal{N}_G$ is a divisor graph by Theorem~\ref{cp-div}.
In the following, we first show that $V(\mathcal{N}_G)\ne  \{p,q,r,s,pq,pr,rs,qs\}$ and $V(\mathcal{N}_G)\ne \{p,q,r,s,pq,pr,ps\}$.

\begin{figure}[hptb]
  \centering
  \includegraphics[width=9cm]{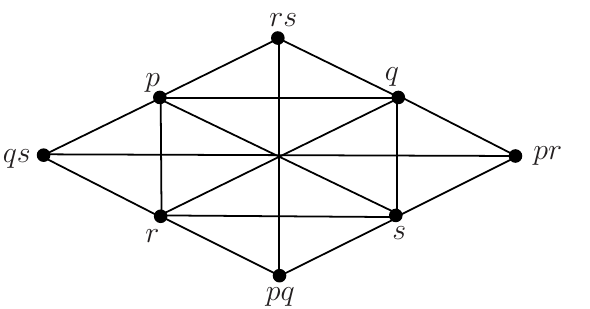}\\
  \caption{$\mathcal{N}_G$ with $V(\mathcal{N}_G)=\{p,q,r,s,pq,pr,rs,qs\}$}\label{t-6}
\end{figure}

Suppose for a contradiction, that $V(\mathcal{N}_G)=\{p,q,r,s,pq,pr,rs,qs\}$. Note that $\mathcal{N}_G$ is a divisor graph. Then
$\mathcal{N}_G$ is shown in Figure~\ref{t-6}. Thus, by Theorem~\ref{dg-o}, $\mathcal{N}_G$ has a transitive orientation, say $\mathcal{T}$. We consider the following two cases.

\smallskip
{\bf Case 1.} $(q,rs)\in E(\mathcal{T})$, which is a directed edge from $q$ to $rs$ in $\mathcal{T}$.
\smallskip

It follows that $(q,pr)\in E(\mathcal{T})$, since $\{pr,rs\}\notin E(\mathcal{N}_G)$. Note that $\{pr,p\}$ is not an edge of $\mathcal{N}_G$. It must be $(q,p)\in E(\mathcal{T})$, and so
\begin{equation}\label{eq-03}
(qs,p)\in E(\mathcal{T}).
\end{equation}
On the other hand, since $\{pq,p\},\{pq,q\}\notin E(\mathcal{N}_G)$, we must have that $(pq,rs)$ belongs to $E(\mathcal{T})$, this implies
\begin{equation}\label{eq-04}
(p,rs)\in E(\mathcal{T}).
\end{equation}
Combining \eqref{eq-03} and \eqref{eq-04}, we obtain $(qs,rs)\in E(\mathcal{T})$, and so $\{qs,rs\}\in E(\mathcal{N}_G)$, which is impossible.

\smallskip
{\bf Case 2.} $(rs,q)\in E(\mathcal{T})$, which is a directed edge from $rs$ to $q$ in $\mathcal{T}$.
\smallskip

Then $(pr,q)\in E(\mathcal{T})$.
Since $(pr,p)\notin E(\mathcal{T})$, it follows that
\begin{equation}\label{eq-05}
(p,q)\in E(\mathcal{T}).
\end{equation}
Also, we deduce that $(rs,pq)\in E(\mathcal{T})$, as $(pq,q)\notin E(\mathcal{T})$. This implies that $(rs,p)\in E(\mathcal{T})$, since $(p,pq)\notin E(\mathcal{T})$. Note that $(rs,qs)\notin E(\mathcal{T})$.
It follows that
\begin{equation}\label{eq-06}
(qs,p)\in E(\mathcal{T}).
\end{equation}
Combining \eqref{eq-05} and \eqref{eq-06}, we have $(qs,q)\in E(\mathcal{T})$, which is impossible as $\{qs,q\}\notin E(\mathcal{N}_G)$.

\begin{figure}[hptb]
  \centering
  \includegraphics[width=5cm]{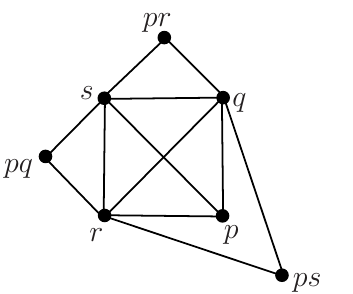}\\
  \caption{$\mathcal{N}_G$ with $V(\mathcal{N}_G)=\{p,q,r,s,pq,pr,ps\}$}\label{t-7}
\end{figure}

We next prove that $V(\mathcal{N}_G)\ne \{p,q,r,s,pq,pr,ps\}$. Suppose otherwise, $V(\mathcal{N}_G)=\{p,q,r,s,pq,pr,ps\}$. Then
$\mathcal{N}_G$ is shown in Figure~\ref{t-7}. In view of Theorem~\ref{dg-o}, we see that $\mathcal{N}_G$ has a transitive orientation, say $\mathcal{G}$. We now consider the following two cases.

\smallskip
{\bf Case 1.} $(s,pr)\in E(\mathcal{G})$, which is a directed edge from $s$ to $pr$ in $\mathcal{G}$.
\smallskip

Then we must have that both $(s,pq)$ and $(s,r)$ belong to  $E(\mathcal{G})$, and so $(s,q)\in E(\mathcal{G})$. This will forces $(ps,q)\in E(\mathcal{G})$, which yields
\begin{equation}\label{eq-07}
(pr,q)\in E(\mathcal{G}).
\end{equation}
On the other hand, since $(s,r)\in E(\mathcal{G})$, it must be $(ps,r)\in E(\mathcal{G})$. This implies $(pq,r)\in E(\mathcal{G})$, and as a result, we have
\begin{equation}\label{eq-08}
(q,r)\in E(\mathcal{G}).
\end{equation}
Combining \eqref{eq-07} and \eqref{eq-08}, we then have $(pr,r)\in E(\mathcal{G})$, a contradiction as $\{pr,r\}\notin E(\mathcal{N}_G)$.

\smallskip
{\bf Case 2.} $(pr,s)\in E(\mathcal{G})$, which is a directed edge from $pr$ to $s$ in $\mathcal{G}$.
\smallskip

It follows that both $(pq,s)$ and $(r,s)$ belong to $E(\mathcal{G})$, and so $(q,s)\in E(\mathcal{G})$. This implies $(q,ps)\in E(\mathcal{G})$, and as a consequence, we obtain
\begin{equation}\label{eq-09}
(q,pr)\in E(\mathcal{G}).
\end{equation}
Also, since $(r,s)\in E(\mathcal{G})$, we have $(r,ps)\in E(\mathcal{G})$. Thus, $(r,pq)\in E(\mathcal{G})$, and then we must have that
\begin{equation}\label{eq-10}
(r,q)\in E(\mathcal{G}).
\end{equation}
Combining \eqref{eq-07} and \eqref{eq-08}, we get a contradiction as $\{r,pr\}\notin E(\mathcal{N}_G)$.

We conclude that $V(\mathcal{N}_G)$ is neither
$\{p,q,r,s,pq,pr,rs,qs\}$ nor $\{p,q,r,s,pq,pr,ps\}$.
As a result, we can choose at most three elements of
$\{pq,pr,qr,ps,qs,rs\}$ belonging to $V(\mathcal{N}_G)$. Furthermore, if we choose precisely three numbers in $\{pq,pr,qr,ps,qs,rs\}$, then these  three numbers can not have a common factor which is a prime. Then
combining Lemma~\ref{teshu}, \eqref{eq-02} and \eqref{eq-0n2}, and without loss of generality, we deduce that \eqref{eq-new2} holds, as desired.
$\qed$

By the proof of Theorem~\ref{4-p-ds} and Lemma~\ref{indsub}, we have the following corollary.

\begin{cor}\label{4-ps}
Given a group $G$, $\Gamma(G)$ is not a divisor graph if
one of the following holds:
\begin{itemize}
  \item[{\rm (a)}] $\{p,q,r,s,pq,pr,rs,qs\}\subseteq \pi_e(G)$;
  \item[{\rm (b)}] $\{p,q,r,s,pq,pr,ps\}\subseteq \pi_e(G)$,
\end{itemize}
where $p,q,r,s$ are pairwise distinct primes.
\end{cor}

\section{Dihedral and Generalized Quaternion Groups}\label{}

Corollary~\ref{nil-g} classifies nilpotent groups $G$ such that
$\Gamma(G)$ is a divisor graph.
In this section, we will classify two families of non-nilpotent groups whose coprime graph is a divisor graph.

Given a positive integer $n$ at least $3$, the {\em dihedral group} of order $2n$, denoted by $D_{2n}$, is defined as follows:
\begin{equation}\label{d2n}
D_{2n}=\langle a,b: b^2=a^n=e, a^{-1}=bab^{-1}\rangle.
\end{equation}
It is easy to see that,
for any $1\le i \le n$, $a^ib$ is an involution, and
\begin{equation}\label{d2n-1}
D_{2n}=\{b,ab,a^2b,\ldots,a^{n-1}b\} \cup \langle a\rangle.
\end{equation}
Recall that $D_{2n}$ is non-nilpotent if and only if $n$ is not a power of $2$.

\begin{thm}
Let $D_{2n}$ be the dihedral group as presented in \eqref{d2n}. Then
$\Gamma(D_{2n})$ is a divisor graph if and only if
$
n=p^tq^m,
$
where $p,q$ are primes and $t,m$ are positive integers.
\end{thm}
\proof
The necessity follows trivially from Corollary~\ref{3-n-pr}.
In the following, suppose that $n=p^tq^m$ for some primes $p,q$ and positive integers $t,m$. Then $|\pi(D_{2n})|\le 3$. Particularly, if $|\pi(D_{2n})|=3$, then $n$ is odd, and so $2p\notin \pi_e(D_{2n})$ by \eqref{d2n-1}, which implies that $\Gamma(D_{2n})$ is a divisor graph from Theorem~\ref{3-p-ds}.
Otherwise, $|\pi(D_{2n})|\le 2$, then it follows from Theorem~\ref{2-p-ds} that $\Gamma(D_{2n})$ is also a divisor graph, as desired.
$\qed$

Given a positive integer $t$, Johnson \cite{Jon} introduced the generalized quaternion group with order $4t$, which is denoted by $Q_{4t}$ and has the following presentation:
\begin{equation}\label{q4t}
Q_{4t}=\langle x,y: x^{2t}=y^4=e, x^t=y^2, y^{-1}xy=x^{-1}\rangle.
\end{equation}
Note that in the literature, $Q_{4t}$ is also called a {\em dicyclic group} and 
$Q_{8}$ is the usual quaternion group of order $8$. Moreover,
the unique involution of $Q_{4t}$ is $y^2=x^t$, and $o(x^iy)=4$ for any $0\le i \le 2t-1$. Note that
\begin{equation}\label{q4t-1}
Q_{4t}=\{x^iy: 0\le i \le 2t-1\} \cup \langle x\rangle.
\end{equation}
Recall that $Q_{4t}$ is non-nilpotent if and only if $t$ is not a power of $2$.

\begin{thm}
Let $Q_{4t}$ be the generalized quaternion group as presented in \eqref{q4t}. Then
$\Gamma(Q_{4t})$ is a divisor graph if and only if
$t=2^mp^n$, where $p$ is a prime and $m,n$ are positive integers.
\end{thm}
\proof
Note that $2t\in \pi_e(Q_{4t})$. Thus, by Corollary~\ref{3-n-pr},
the necessity is valid.
For the converse, let $t=2^mp^n$ for some prime $p$ and positive integers $m,n$. It follows from \eqref{q4t-1} that $|\pi(Q_{4t})|\le 2$, and so
by Theorem~\ref{2-p-ds}, we have that $\Gamma(Q_{4t})$ is also a divisor graph, as desired.
$\qed$

\section{Symmetric Groups and Alternating Groups}\label{}

In this section, we will classify the symmetric groups and alternating groups whose coprime graph is a divisor graph.

Recall that for every $n\ge 4$, $\mathbf{A}_n$ is  non-nilpotent, and if $n\ge 3$, then $\mathbf{S}_n$ is non-nilpotent.
As we all know, for any $n\ge 5$, $\mathbf{A}_n$ is simple.
Symmetric group is very important in diverse areas of mathematics, such as combinatorics and the representation theory of Lie groups.

\begin{thm}\label{sym-g}
The following hold:
\begin{itemize}
  \item[\rm (1)] $\Gamma(\mathbf{S}_n)$ is a divisor graph if and only if $n\le 7$;
  \item[\rm (2)] $\Gamma(\mathbf{A}_n)$ is a divisor graph if and only if $n\le 8$.
\end{itemize}
\end{thm}
\proof
(1) We first prove that $\Gamma(\mathbf{S}_8)$ is not a divisor graph.
Note that
$$
\{(1,2)(3,4,5),(1,2)(3,4,5,6,7),(1,2,3)(4,5,6,7,8)\}\subseteq\mathbf{S}_8.
$$
As a result, we have that $\{6,10,15\}\subseteq \pi_e(\mathbf{S}_8)$, and so $\Gamma(\mathbf{S}_8)$ is not a divisor graph by Lemma~\ref{teshu}, as desired. Since $\mathbf{S}_{k}$ is a subgroup $\mathbf{S}_{l}$ for any two integers $k,l$ with $k\le l$. Now in view of Observation~\ref{obs-1}, it suffices to prove that $\Gamma(\mathbf{S}_7)$ is a divisor graph.
Note that $\pi_e(\mathbf{S}_7)=\{2, 3, 4, 5, 6, 7, 10, 12\}$. Therefore, it follows that $V(\mathcal{N}_{\mathbf{S}_7})=\{2,3,5,7,6,10\}$, and so
$\Gamma(\mathbf{S}_7)$ is a divisor graph by Theorem~\ref{4-p-ds}, as desired.

(2) We first prove that $\Gamma(\mathbf{A}_9)$ is not a divisor graph.
Note that
$$
\{(1,2)(3,4,5)(6,7),(1,2)(3,4,5,6,7)(8,9),(1,2,3)(4,5,6,7,8)\}\subseteq \mathbf{A}_9.
$$
Thus, $\{6,10,15\}\subseteq \pi_e(\mathbf{A}_9)$ and so, Lemma~\ref{teshu} implies that $\Gamma(\mathbf{A}_9)$ is not a divisor graph, as desired. It follows from Observation~\ref{obs-1} that,
it suffices to prove that $\Gamma(\mathbf{A}_8)$ is a divisor graph.
Note that $\pi_e(\mathbf{A}_8)=\{2, 3, 4, 5, 6, 7, 15\}$.
As a result, $|\pi(\mathbf{A}_8)|=4$ and $V(\mathcal{N}_{\mathbf{A}_8})=\{2,3,5,7,6,15\}$.
It follows from Theorem~\ref{4-p-ds} that $\Gamma(\mathbf{A}_8)$ is a divisor graph, as desired.
$\qed$

\section{Direct Products}\label{}

For two non-trivial groups $H$ and $K$, for which direct product $H\times K$ is the coprime graph a divisor graph?
In this section, we describe the direct product of two non-trivial groups whose coprime graph is a divisor graph. Our main result is the following theorem.

\begin{thm}\label{zhiji}
Given  two non-trivial groups $H$ and $K$, $\Gamma(H\times K)$ is a divisor graph if and only if one of the following holds:
\begin{itemize}
\item[{\rm (a)}] $\pi(H)=\pi(K)=\{p,q\}$ for distinct primes $p$ and $q$;

\item[{\rm (b)}] $\pi(H)=\{p\}$, $\pi(K)\subseteq\{p,q,r\}$ and $qr\notin \pi_e(K)$, where $p,q,r$ are pairwise distinct primes.
\end{itemize}
\end{thm}
\proof
Let $G=H\times K$.
If (a) occurs, then $|\pi(G)|=2$, and so $\Gamma(G)$ is a divisor graph by Theorem~\ref{2-p-ds}. Now suppose that (b) occurs. Then $|\pi(G)|\le 3$.
If one of $q$ and $r$ does not belong to $\pi(K)$, then $|\pi(G)|\le 2$, it follows from Theorem~\ref{2-p-ds} that $\Gamma(G)$ is a divisor graph. Thus, in the following, we may assume that $\{q,r\}\subseteq\pi(K)$.
It follows that $\pi(G)=\{p,q,r\}$. Note that the fact that, if $G$ has an element $x$ with $\pi(x)=\{q,r\}$, then $x\in K$. Since $qr\notin \pi_e(K)$, we have that $qr\notin \pi_e(G)$.
Now Theorem~\ref{3-p-ds} implies that $\Gamma(G)$ is a divisor graph.

For the converse, suppose that $\Gamma(G)$ is a divisor graph.
In the following, we consider two cases.

\smallskip
{\bf Case 1.} Neither $H$ nor $K$ is a $p$-group, where $p$ is a prime.
\smallskip

Suppose for a contradiction that $|\pi(H)|\ge 3$.
Let $\{p,q,r\}\subseteq\pi(H)$, where $p,q,r$ are pairwise distinct primes. If there exists a prime $s\in \pi(K)$ such that $s\notin \{p,q,r\}$, then $\{ps,qs,rs\}\subseteq \pi_e(G)$, and so from Corollary~\ref{4-ps}, it follows that $\Gamma(G)$ is not a divisor graph, a contradiction.
As a result, we have $\pi(K)\subseteq \{p,q,r\}$. Note that $K$ is not a $p$-group. Without loss of generality, let $\{p,q\}\subseteq \pi(K)$. Then $\{pq,pr,qr\}\subseteq \pi_e(G)$, and thus, $\Gamma(G)$ is not a divisor graph by Lemma~\ref{teshu}, a contradiction. Consequently, we conclude that $|\pi(H)|=2$, and similarly, we also have $|\pi(K)|=2$.

Let now $\pi(H)=\{p,q\}$ with distinct primes $p,q$. If $\pi(H)\cap \pi(K)=\emptyset$, say $\pi(K)=\{r,s\}$ with distinct primes $r,s$, then $\{pr,ps,qr,qs\}\in \pi_e(G)$, and so $\Gamma(G)$ is not a divisor graph by Corollary~\ref{4-ps}, a contradiction. If $|\pi(H)\cap \pi(K)|=1$, say $\pi(K)=\{p,r\}$ for prime $r\notin \{p,q\}$, then $\{pq,pr,qr\}\subseteq \pi_e(G)$, it follows from Lemma~\ref{teshu} that $\Gamma(G)$ is not a divisor graph, a contradiction. As a consequence, in this case, it must be that $\pi(H)=\pi(K)=\{p,q\}$, this implies that (a) holds.

\smallskip
{\bf Case 2.} One of $H$ and $K$ is a $p$-group, where $p$ is a prime.
\smallskip

Without loss of generality, assume that $\pi(H)=\{p\}$. If $\{q,r,s\}\subseteq\pi(K)$ for pairwise distinct primes $q,r,s$, then $\{pq,pr,ps\}\subseteq \pi_e(G)$, a contradiction as Corollary~\ref{4-ps}.
It follows that $\pi(K)\subseteq\{p,q,r\}$ with pairwise distinct primes $p,q,r$. If $qr\in \pi_e(K)$, then it must be that $\{pq,pr,qr\}\subseteq \pi_e(G)$, which is impossible by  Lemma~\ref{teshu}. Therefore, in this case, (b) holds.
$\qed$

By Theorem~\ref{zhiji}, we have the following examples.

\begin{example}\label{ex-zj}
$\Gamma(G)$ is not a divisor graph if $G$ is isomorphic to one of the following:
$$
\mathbf{S}_{3}\times \mathbf{A}_{5},D_{6}\times D_{10}, \mathbf{A}_{5}\times D_{10},\mathbf{S}_5\times L_3(2),\mathbb{Z}_2 \times Sz(8),\mathbb{Z}_3\times G_2(3),\mathbb{Z}_2\times M_{11}.
$$
\end{example}

\section{Sporadic Simple Groups}\label{}

In this section, we will classify all sporadic simple groups whose coprime graph is a divisor graph. Our main result is the following theorem.

\begin{thm}\label{ssg}
Suppose that $G$ is a sporadic simple group. Then $\Gamma(G)$ is a divisor graph if and only if $G$ is isomorphic to one of the following Mathieu groups: $$M_{11},~~M_{12},~~M_{22},~~M_{23}.$$
\end{thm}

We next show two lemmas before giving the proof of Theorem~\ref{ssg}.

\begin{lemma}\label{dg-le-1}
For $n\ge 3$,
suppose that $V(\mathcal{N}_G)=\{p_1,p_2,\ldots,p_n,p_lp_t\}$, where $p_i$ is a prime for every $1\le i \le n$ and $1\le l < t \le n$.
Then $\Gamma(G)$ is a divisor graph.
\end{lemma}
\proof
Without loss of generality, we may assume that $l=n-2$ and $t=n-1$.
Let $\Delta$ be the subgraph of $\mathcal{N}_G$ induced by $\{p_1,p_2,\ldots,p_n\}$. Then $\Delta$ is complete, and so $\Delta$ has a transitive orientation by Lemma~\ref{complete} and Theorem~\ref{dg-o}.
Now we define a transitive orientation $\mathcal{T}_{\Delta}$ of $\Delta$ as follows:
\begin{equation}\label{eq-dg-1}
(p_i,p_j)\in E(\mathcal{T}_{\Delta}) \Leftrightarrow 1\le i<j \le n.
\end{equation}
Namely, there is a directed edge from $p_i$ to $p_j$ in $\mathcal{T}_{\Delta}$ if and only if $i<j$.
Note that $\mathcal{N}_G$ is the graph obtained by adding the edges
\begin{equation}\label{eq-dg-2}
\{\{p_i,p_lp_t\}: i=1,2,\dots,n-3,n\}
\end{equation}
and the vertex $p_lp_t$ into $\Delta$. Now we give an orientation $\mathcal{T}_{\mathcal{N}_G}$ for $\mathcal{N}_G$ as follows:
\begin{itemize}
  \item $E(\mathcal{T}_{\Delta})\subseteq E(\mathcal{T}_{\mathcal{N}_G})$;
  \item $(p_i,p_lp_t)\in E(\mathcal{T}_{\mathcal{N}_G})$, $i=1,2,\ldots,n-3$;
  \item $(p_lp_t,p_n)\in E(\mathcal{T}_{\mathcal{N}_G})$.
\end{itemize}
In the following, we prove that $\mathcal{T}_{\mathcal{N}_G}$ is a transitive orientation. Note that $p_n$ is a receiver in  $\mathcal{T}_{\mathcal{N}_G}$ and $\mathcal{T}_{\Delta}$ is a transitive orientation of $\Delta$. By \eqref{eq-dg-2},
it suffices to show that if $(a,b),(b,p_lp_t)\in \mathcal{T}_{\mathcal{N}_G}$ for distinct $a,b\in V(\Delta)$, then $(a,p_lp_t)\in \mathcal{T}_{\mathcal{N}_G}$.
Now assume that $(a,b),(b,p_lp_t)\in \mathcal{T}_{\mathcal{N}_G}$ for distinct $a,b\in V(\Delta)$. By the definition of  $\mathcal{T}_{\mathcal{N}_G}$, we have that $b=p_k$ for some $1\le k \le n-3$. Then it follows from \eqref{eq-dg-1} that $k\ge 2$ and
$a=p_r$ for some $1\le r < k$. This forces that $(p_r,p_lp_t)\in \mathcal{T}_{\mathcal{N}_G}$, and so $\mathcal{T}_{\mathcal{N}_G}$ is a transitive orientation. According to Theorem~\ref{cp-div}, we have that
$\Gamma(G)$ is a divisor graph.
$\qed$

\begin{figure}[hptb]
  \centering
  \includegraphics[width=8cm]{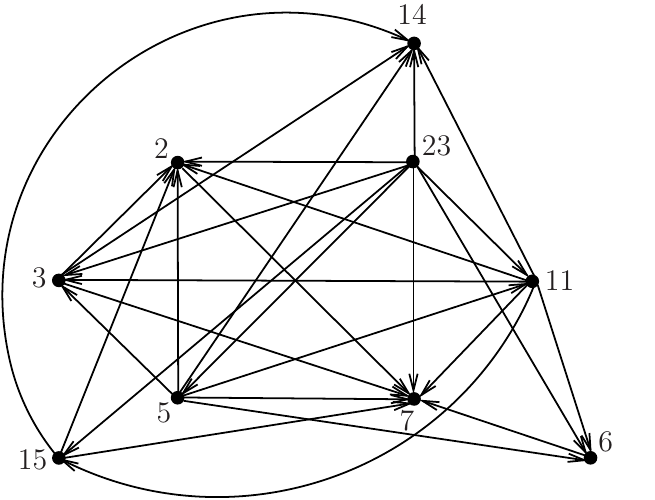}\\
  \caption{A transitive orientation for $\mathcal{N}_{M_{23}}$}\label{t-8}
\end{figure}

\begin{lemma}\label{M23}
$\Gamma(M_{23})$ is a divisor graph.
\end{lemma}
\proof
It is well known that
$\pi_e(M_{23})=\{2, 3, 4, 5, 6, 7, 8, 11,14,15,23\}$,
and so
$$
\mathcal{N}_{M_{23}}=\{2,3,5,7,11,23,6,14,15\}.
$$
It is easy to verify that Figure~\ref{t-8} is a transitive orientation of $\mathcal{N}_{M_{23}}$, and so $\Gamma(M_{23})$ is a divisor graph by Theorem~\ref{dg-o}.
$\qed$

We are now ready to prove Theorem~\ref{ssg}.
\medskip

\noindent {\em Proof of Theorem~{\rm\ref{ssg}}.}
As we all know, there exist exactly $26$ sporadic simple groups.
We first consider the five Mathieu groups.
For $M_{11}$ and $M_{12}$, we have that
$$
\pi_e(M_{11})=\{2, 3, 4, 5, 6, 8, 11 \},
~~\pi_e(M_{12})=\{2, 3, 4, 5, 6, 8, 10, 11 \}.
$$
Then $|\pi(M_{11})|=|\pi(M_{12})|=4$, and by Theorem~\ref{4-p-ds}, we see that both $\Gamma(M_{12})$ and $\Gamma(M_{22})$ are divisor graphs.
For $M_{22}$, we have $\pi_e(M_{22})=\{1, 2, 3, 4, 5, 6, 7, 8, 11\}$. As a result, $\mathcal{N}_{M_{22}}=\{2,3,5,7,11,6\}$, and so by Lemma~\ref{dg-le-1}, we see that $\Gamma(M_{22})$ is a divisor graph.
For $M_{23}$, from Lemma~\ref{M23}, it follows that $\Gamma(M_{23})$ is a divisor graph. Also,
note that the fact that $\{6,10,15\}\in \pi_e(M_{24})$. We deduce that  $\Gamma(M_{24})$ is not a divisor graph by Lemma~\ref{teshu}.

Then by \cite{CC85}, it follows that Janko group $J_1$, Janko group $J_2$, Janko group $J_4$, Held group $He$, Harada-Norton group $HN$, Thompson group $Th$, Baby Monster group $B$, Monster group $M$, O'Nan group $O'N$, Lyons group $Ly$, Rudvalis group $Ru$, Suzuki group $Suz$, Fischer group $Fi_{22}$,
and Higman-Sims group $HS$
contain
$D_6\times D_{10}$, $\mathbf{A}_5\times D_{10}$, $M_{24}$, $\mathbf{S}_4\times L_3(2)$, $\mathbf{A}_{12}$, $\mathbb{Z}_3\times G_2(3)$, $Th$, $\mathbf{A}_{12}$, $J_1$, $\mathbb{Z}_2\times M_{11}$, $\mathbb{Z}_2\times \mathbb{Z}_2 \times Sz(8)$, $\mathbf{S}_{3}\times \mathbf{A}_{5}$, $\mathbf{S}_{10}$, and $\mathbf{S}_8$
as subgroups, respectively.
By Example~\ref{ex-zj} and Theorem~\ref{sym-g},
we have that every of the above simple groups has a non-divisor coprime graph.

Finally, we handle the remaining $7$ sporadic simple groups.
Note first that
$$\{6,10,15\}\in \pi_e(Mcl)\cap \pi_e(J_3) \cap \pi_e(Fi_{23})\cap \pi_e(Fi_{24}').$$
Therefore, every of $\Gamma(Mcl)$, $\Gamma(J_3)$, $Fi_{23}$ and $Fi_{24}'$ has a non-divisor coprime graph by Lemma~\ref{teshu}.
Then, by \cite{CC85}, it is easy to see that
for every $1\le i \le 3$,
Conway group $Co_i$ has a subgroup isomorphic to $McL$ .
As a result, by Observation~\ref{obs-1}, we have that $\Gamma(Co_i)$ is not a divisor graph for each $1\le i \le 3$.
$\qed$

\section{Conclusions}

This paper mainly studies the coprime graph of a finite group  $G$ and is motivated by the following question: which coprime graphs are divisor graphs? 
More precisely, 
this paper first characterizes which generalized lexicographic products are divisor graphs. 
As applications, it shows that power graphs, reduced power graphs and order graphs are all divisor graphs, which also implies the main result in \cite{Ta}. 
Then, this paper shows that the coprime graph of a group is a generalized 
lexicographic product, and characterizes which coprime graphs are divisor graphs. Finally, the paper classifies the finite groups $G$ having at most four prime divisors, whose coprime graphs are divisor graphs, and also classifies the finite groups $G$ whose coprime graphs are divisor graphs provided that $G$ is a nilpotent group, a dihedral group, a generalized quaternion group, a symmetric group, an alternating group, a direct product of two non-trivial groups, and a sporadic simple group.

In the future, we will further study finite groups whose coprime graph is a divisor graph. For example, we will attempt to classify all finite simple groups whose coprime graph is a divisor graph.

\section*{Acknowledgements}
We are grateful to the referees for detailed comments and thoughtful suggestions that have helped improve this paper substantially. 

This research was supported by the National Natural Science Foundation of China (Grant No. 12326333) and the
Shaanxi Fundamental Science Research Project for Mathematics and Physics (Grant No. 22JSQ024).

\section*{Use of AI tools declaration}
The authors declare they have not used Artificial Intelligence (AI) tools in the creation of this article.

\section*{Conflict of interest}

The authors declare there is no conflicts of interest.

\end{document}